\documentclass[twoside,11pt]{article}
\usepackage{pgm}

\usepackage{amsmath, amssymb}

\usepackage[utf8]{inputenc}
\inputencoding{utf8}

\usepackage{hyperref,color,soul,booktabs}
\setulcolor{blue}
\newcommand{\BibTeX}{\textsc{B\kern-0.1emi\kern-0.017emb}\kern-0.15em\TeX}

\newcommand{\sbsbsection}[1]{\noindent \textbf{#1.}  }

\usepackage[usenames, dvipsnames]{xcolor}

\ShortHeadings{The covariance of causal effect estimators for binary v-structures}{Kuipers and Moffa}

\begin{document}

\title{The covariance of causal effect estimators for binary v-structures}

\author{\Name{Jack Kuipers} \Email{jack.kuipers@bsse.ethz.ch} \\
\addr D-BSSE, ETH Zurich, Schanzenstrasse 44, 4056 Basel, Switzerland
\and
\Name{Giusi Moffa} \Email{giusi.moffa@unibas.ch} \\
\addr Department of Mathematics and Computer Science, University of Basel, Basel, Switzerland}

\maketitle

\begin{abstract}
Previously [Journal of Causal Inference, 10, 90--105 (2022)], we computed the variance of two estimators of causal effects for a v-structure of binary variables. Here we show that a linear combination of these estimators has lower variance than either. Furthermore, we show that this holds also when the treatment variable is block randomised with a predefined number receiving treatment, with analogous results to when it is sampled randomly.
\end{abstract}

\section{Recap}

\begin{figure}[t]
 \centering
 \includegraphics[width=0.4\textwidth]{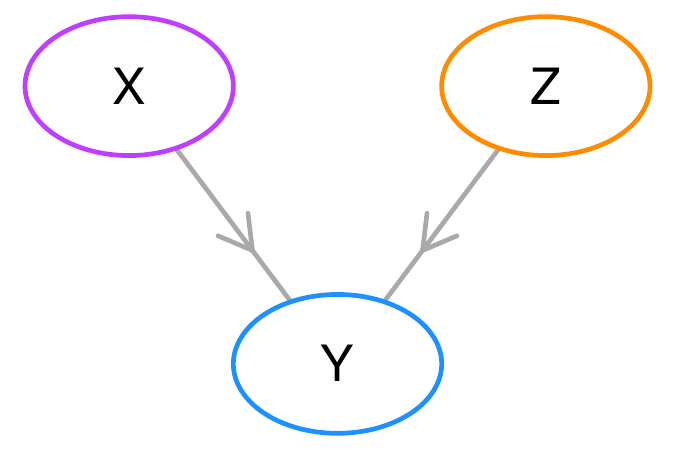}
\caption{A v-structure on 3 nodes.}
\label{fig:vstruct}
\end{figure}

In \cite{km22} we considered a DAG with 3 nodes of binary variables organised in a v-structure with the outcome $Y$ of interest as a collider with parents $Z$ and $X$ (Figure \ref{fig:vstruct}) where we want to estimate the effect of $X$ on $Y$. We recall that the DAG in Figure \ref{fig:vstruct} has the following probability tables:
\begin{align}
p(X = 1) = p_X \, , & & p(Y = 1 \mid X = 0, Z = 0) &= p_{Y, 0}  \, , & & p(Y = 1 \mid X = 1, Z = 0) = p_{Y, 2} \nonumber \\
p(Z = 1) = p_Z \, , & & p(Y = 1 \mid X = 0, Z = 1) &= p_{Y, 1}  \, , & & p(Y = 1 \mid X = 1, Z = 1) = p_{Y, 3}
\end{align}
which is equivalent to multinomial sampling with probabilities
\begin{equation}
\begin{array}{c|c|c|lcc|c|c|l}
X & Z & Y & p & & X & Z & Y & p \\
\cline{1-4}
\cline{6-9}
0 & 0 & 0 & p_0 = (1-p_X)(1-p_Z)(1-p_{Y,0}) & & 1 & 0 & 0 & p_4 = p_X(1-p_Z)(1-p_{Y,2}) \\
0 & 0 & 1 & p_1 = (1-p_X)(1-p_Z)p_{Y,0} & & 1 & 0 & 1 & p_5 = p_X(1-p_Z)p_{Y,2} \\
0 & 1 & 0 & p_2 = (1-p_X)p_Z(1-p_{Y,1}) & & 1 & 1 & 0 & p_6 = p_Xp_Z(1-p_{Y,3}) \\
0 & 1 & 1 & p_3 = (1-p_X)p_Zp_{Y,1} & & 1 & 1 & 1 & p_7 = p_Xp_Zp_{Y,3}
\end{array}
\end{equation}

If we represent with $N_i$ the number of sampled binary vectors indexed by $i = 4X+2Z+Y$, with a total of $N$, the causal estimator from raw conditionals is 
\begin{align}
R = R_1 - R_0\, , & & R_1 = \frac{N_5 + N_7}{N_4 + N_5 + N_6 + N_7}\, , & & R_0 = \frac{N_1 + N_3}{N_0 + N_1 + N_2 + N_3}
\end{align}
while the estimator using marginalisation is
\begin{align}
M = M_1 - M_0 \, , & &
M_1 = M_{11} + M_{10} \, , & &  M_0 = M_{01} + M_{00}
\end{align}
with the terms separated for later ease
\begin{align}
M_{11} = \frac{N_7}{(N_6+N_7)}\frac{(N_2 + N_3 + N_6 + N_7)}{N} \, , & & M_{01} = \frac{N_3}{(N_2+N_3)}\frac{(N_2 + N_3 + N_6 + N_7)}{N} \nonumber \\
M_{10} = \frac{N_5}{(N_4+N_5)}\frac{(N_0 + N_1 + N_4 + N_5)}{N} \, , & & M_{00} = \frac{N_1}{(N_0+N_1)}\frac{(N_0 + N_1 + N_4 + N_5)}{N}
\end{align}

\section{A combined estimator}

We previously computed the variances of these two estimators \citep{km22}, but since both $R$ and $M$ are estimators of the effect of $X$ on $Y$, so too is any linear combination of them
\begin{align}
P(\alpha) = \alpha R + (1-\alpha)M
\end{align}
We could even have different linear combinations for the different parts of the estimators
\begin{align}
Q = \alpha_1 R_1 + (1-\alpha_1)M_1 -  \alpha_0 R_0 - (1-\alpha_0)M_0  
\end{align}
though we will stick with the simpler combination here.

The variance of $P$ is given by
\begin{align}
V[P(\alpha)] = \alpha^2V[R] + (1-\alpha)^2V[M] + 2\alpha(1-\alpha)C[R,M]
\end{align}
in terms of the covariance of $R$ and $M$. Notably, if the covariance is smaller than both variances, then the variance of $P$ is minimised at
\begin{align}
\alpha^* = \frac{V[M] - C[R,M]}{(V[R] - C[R, M]) + (V[M] - C[R, M])}
\end{align}
and the variance is lower than taking the either the $R$ or $M$ estimator:
\begin{align}
C[R, M] < \min(V[R], V[M]) \implies  V[P(\alpha^*)] < \min(V[R], V[M])
\end{align}

\section{Computing the covariance}

To compute the covariance we need different generating variables than before \citep{km22} so we introduce the following generating function:
\begin{align}
U_N = & \left\{[p_0s_1t_1+p_2s_2t_2 + (p_1s_1t_1u_1+p_3s_2t_2u_2)v_1]w_1 + \right. \\ \nonumber
& \quad \left. [p_4s_1t_3+p_6s_2t_4 + (p_5s_1t_3u_3+p_7s_2t_4u_4)v_2]w_2\right\}^{N}
\end{align}
Let $G = \{s_1, s_2, t_1, t_2, t_3, t_4, u_1, u_2, u_3, u_4, v_1, v_2, w_1, w_2 \}$ be all generating variables.

Using this formulation
\begin{align} \label{M11R1expected}
E[M_{11}R_{1}] &= \frac{1}{N}\int \mathrm{d}w_2 \int\mathrm{d}t_4 \frac{s_2u_4v_2}{t_4w_2}\frac{\partial^3}{\partial s_2\partial u_4 \partial v_2} U_N \, \Bigg|_{\substack{G_i = 1 \cr \forall i}} \nonumber \\
&= 
\int \mathrm{d}w_2 \frac{v_2}{w_2}\frac{\partial}{\partial v_2} \frac{(p_3v_1w_1 + p_7v_2w_2 + p_2w_1 + p_6w_2)v_2p_7}{(p_6 + p_7v_2)} U_{N-1} \, \Bigg|_{\substack{G_i = 1 \cr \forall i}} \nonumber \\
&= \frac{p_6p_7}{(p_6+p_7)^2}(N-1)p_Zp_X(1-p_X)^{N-2}F\left([1, 1, 2-N], [2, 2], -\frac{p_X}{1-p_X}\right) + \nonumber \\
& \qquad \frac{(p_4 + p_6)p_7}{Np_{X}^2} + E[M_{11}]E[R_1] 
\end{align}
and likewise
\begin{align} \label{M11R0expected}
E[M_{11}R_{0}] &= 
\int \mathrm{d}w_1 \frac{v_1}{w_1}\frac{\partial}{\partial v_1} \frac{(p_3v_1w_1 + p_7v_2w_2 + p_2w_1 + p_6w_2)v_2p_7}{(p_6 + p_7v_2)} U_{N-1} \, \Bigg|_{\substack{G_i = 1 \cr \forall i}} \nonumber \\
&=  \frac{(p_0p_3 - p_1p_2)p_7}{N(1-p_{X})^2(p_6+p_7)} + E[M_{11}]E[R_0] 
\end{align}

These calculations can be directly repeated to compute the covariances between $M_{10}$ and $R$:
\begin{align} \label{M10R1expected}
E[M_{10}R_{1}] &= \frac{1}{N}\int \mathrm{d}w_2 \int\mathrm{d}t_3 \frac{s_1u_3v_2}{t_3w_2}\frac{\partial^3}{\partial s_1\partial u_3 \partial v_2} U_N \, \Bigg|_{\substack{G_i = 1 \cr \forall i}} \nonumber \\
&= 
\int \mathrm{d}w_2 \frac{v_2}{w_2}\frac{\partial}{\partial v_2} \frac{(p_1v_1w_1 + p_5v_2w_2 + p_0w_1 + p_4w_2)v_2p_5}{(p_4 + p_5v_2)} U_{N-1} \, \Bigg|_{\substack{G_i = 1 \cr \forall i}} \nonumber \\
&= \frac{p_4p_5}{(p_4+p_5)^2}(N-1)(1-p_Z)p_X(1-p_X)^{N-1}F\left([1, 1, 2-N], [2, 2], -\frac{p_X}{1-p_X}\right) + \nonumber \\
& \qquad \frac{(p_4 + p_6)p_5}{Np_{X}^2} + E[M_{10}]E[R_1] 
\end{align}
and
\begin{align} \label{M10R0expected}
E[M_{10}R_{0}] &= 
\int \mathrm{d}w_1 \frac{v_1}{w_1}\frac{\partial}{\partial v_1} \frac{(p_1v_1w_1 + p_5v_2w_2 + p_0w_1 + p_4w_2)v_2p_5}{(p_4 + p_5v_2)} U_{N-1} \, \Bigg|_{\substack{G_i = 1 \cr \forall i}} \nonumber \\
&=  -\frac{(p_0p_3 - p_1p_2)p_5}{N(1-p_{X})^2(p_4+p_5)} + E[M_{10}]E[R_0] 
\end{align}

Similarly for the covariances of $M_{00}$ and $R$:
\begin{align} \label{M00R1expected}
E[M_{00}R_{1}] &=  \frac{(p_5p_6 - p_4p_7)p_1}{Np_{X}^2(p_0+p_1)} + E[M_{00}]E[R_1] 
\end{align}
\begin{align} \label{M00R0expected}
E[M_{00}R_{0}] &= \frac{p_0p_1}{(p_0+p_1)^2}(N-1)(1-p_Z)(1-p_X)p_X^{N-1}F\left([1, 1, 2-N], [2, 2], -\frac{1-p_X}{p_X}\right) + \nonumber \\
& \qquad \frac{(p_0 + p_2)p_1}{N(1-p_{X})^2} + E[M_{00}]E[R_0] 
\end{align}
And for the covariances of $M_{01}$ and $R$:
\begin{align} \label{M01R1expected}
E[M_{01}R_{1}] &=  -\frac{(p_5p_6 - p_4p_7)p_3}{Np_{X}^2(p_2+p_3)} + E[M_{01}]E[R_1] 
\end{align}
\begin{align} \label{M01R0expected}
E[M_{01}R_{0}] &= \frac{p_2p_3}{(p_2+p_3)^2}(N-1)p_Z(1-p_X)p_X^{N-1}F\left([1, 1, 2-N], [2, 2], -\frac{1-p_X}{p_X}\right) + \nonumber \\
& \qquad \frac{(p_0 + p_2)p_3}{N(1-p_{X})^2} + E[M_{01}]E[R_0] 
\end{align}

The complete result is then
\begin{align}
& C[R, M] \cdot N = \frac{p_6p_7}{(p_6+p_7)^2}N(N-1)p_Zp_X(1-p_X)^{N-2}F\left([1, 1, 2-N], [2, 2], -\frac{p_X}{1-p_X}\right) \nonumber \\ 
& + \frac{p_4p_5}{(p_4+p_5)^2}N(N-1)(1-p_Z)p_X(1-p_X)^{N-1}F\left([1, 1, 2-N], [2, 2], -\frac{p_X}{1-p_X}\right) \nonumber \\ 
& + \frac{p_2p_3}{(p_2+p_3)^2}N(N-1)p_Z(1-p_X)p_X^{N-1}F\left([1, 1, 2-N], [2, 2], -\frac{1-p_X}{p_X}\right) \nonumber \\ 
& + \frac{p_0p_1}{(p_0+p_1)^2}N(N-1)(1-p_Z)(1-p_X)p_X^{N-1}F\left([1, 1, 2-N], [2, 2], -\frac{1-p_X}{p_X}\right) \nonumber \\ 
& + \frac{(p_4 + p_6)(p_5 + p_7)}{p_{X}^2} + \frac{(p_0 + p_2)(p_1 + p_3)}{(1-p_{X})^2} \nonumber \\
& + 2\left(\frac{p_7}{(p_6+p_7)} - \frac{p_5}{(p_4+p_5)}\right)\left(\frac{p_3}{(p_2+p_3)} - \frac{p_1}{(p_0+p_1)}\right)p_{Z}(1-p_{Z}) 
\end{align}

\subsection{Numerical checks}

We include code to evaluate the analytical results and run simulations at \url{https://github.com/jackkuipers/Vcausal}. For example, for $p_X = \frac{1}{3}, p_Z = \frac{2}{3}, p_{Y,0} = \frac{1}{6}, p_{Y,1} = \frac{1}{2}, p_{Y,2} = \frac{1}{3}, p_{Y,3} = \frac{5}{6}$ and $N=100$, we obtained Monte Carlo estimates of the square root of the covariance of $R$ and $M$ as $0.0915292$ from 40 million repetitions, in line with the analytical result of $0.0915308$. This is lower than the analytical results for $R$ and $M$ of $0.1019324$ and $0.0924017$ respectively (with MC estimates of $0.1019294$ and $0.0924014$). This means that the best combination would have $\alpha^*=0.0737371$, for which we would have a minimum variance of $0.0923378$. This aligns with the result from the simulation of $0.0923373$.

\subsection{Relative difference in variances}

As in \cite{km22}, we use the following parametrisation
\begin{align} \label{intersetting}
p_{Y, 0} &= q_0 - C - Dp_{Z}\, , & p_{Y, 1} &= q_0 + C + D(1-p_{Z})\, , \nonumber \\
p_{Y, 2} &= q_1 - C + Dp_{Z}\, , & p_{Y, 3} &= q_1 + C - D(1-p_{Z}) \, ,
\end{align}
so that the causal effect of $X$ on $Y$ is $q_1 - q_0$. The effect of $Z$ on $Y$ in the probability tables is then $2C+D$ for $X=0$ and $2C-D$ for $X=1$, making $C$ a measure of the causal effect of $Z$ on $Y$ and $D$ a measure of the interaction effect between $X$ and the effect of $Z$ on $Y$.

\begin{figure}
 \centering
 \begin{tabular}{cc}
 \includegraphics[width=0.45\textwidth]{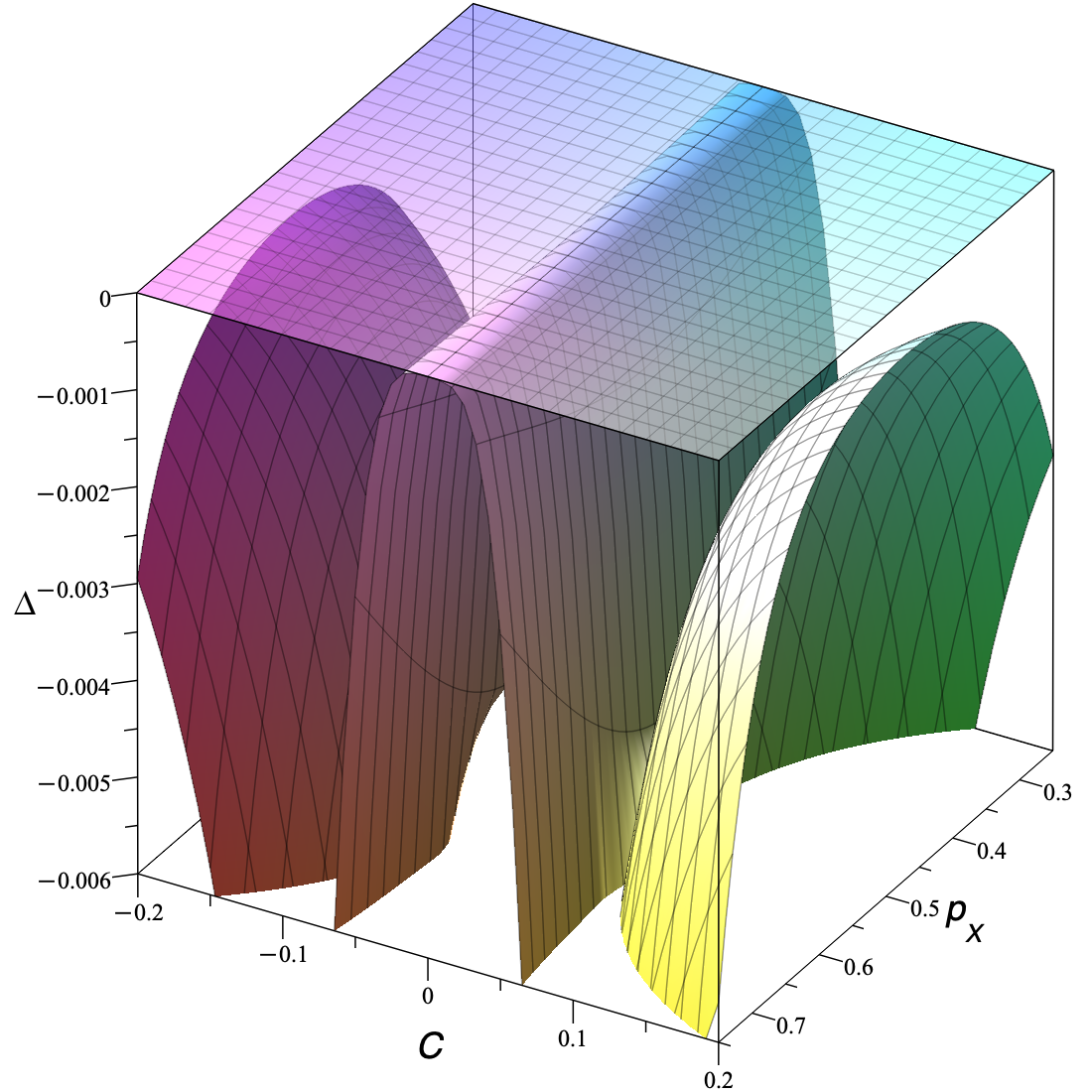} & \includegraphics[width=0.45\textwidth]{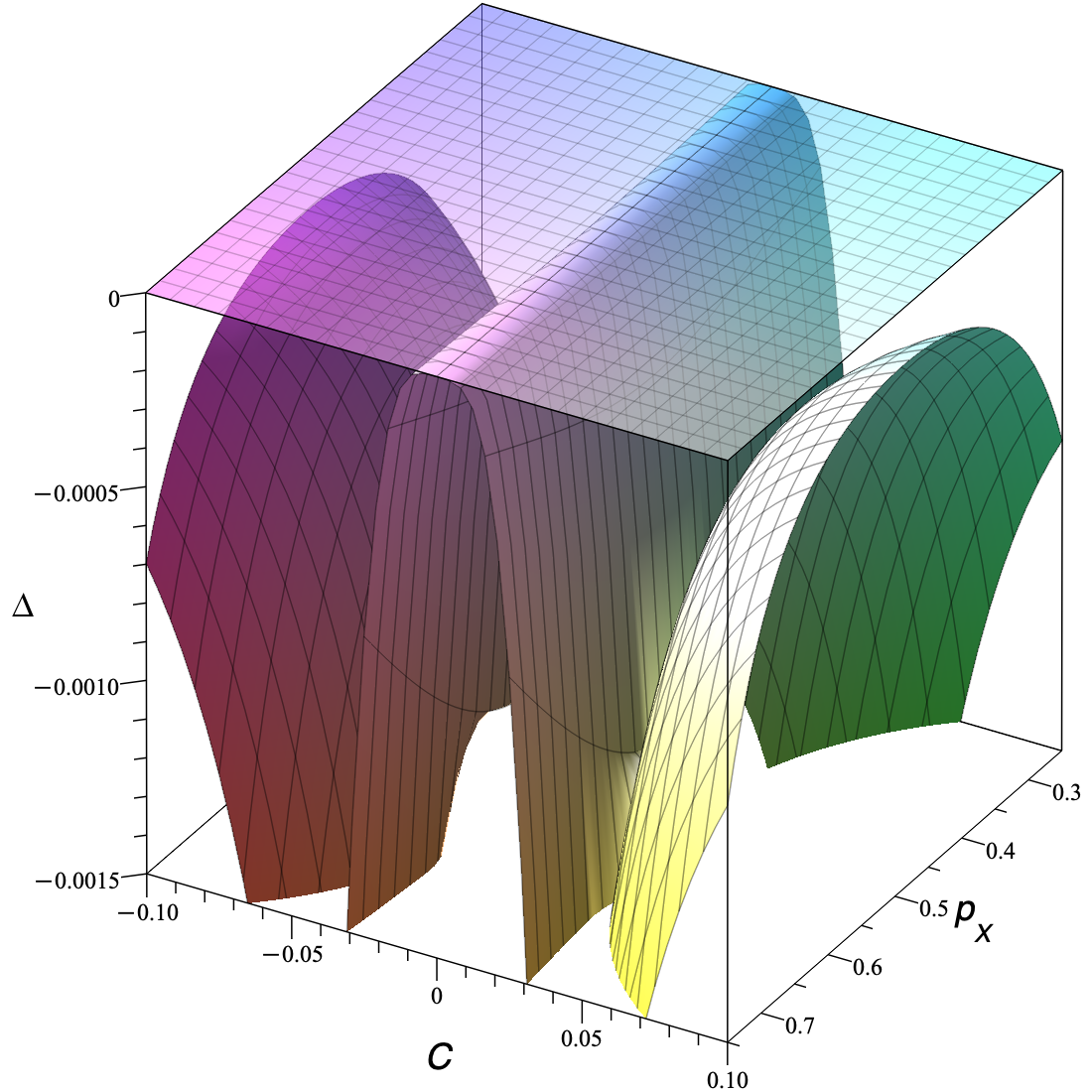} \\
 (a) $N=100$ & (b) $N=400$
 \end{tabular}
\caption{Relative difference in variance of the combined estimator to the better of the other two.}
\label{fig:deltacoplot}
\end{figure}

We plot the difference in variances of the combined estimator to the better of the other two estimators:
\begin{align}
\Delta = \frac{V[P(\alpha^*)]- \min(V[M],V[R])}{\min(V[M],V[R])}
\end{align}
In Figure \ref{fig:deltacoplot} we leave $p_X$ free, set $p_Z=\frac{2}{3}$ and $D=0$, and set $q_0=\frac{1}{3}$, $q_1=\frac{2}{3}$. We plot $\Delta$ for $N=100$ and $N=400$, where in the plot for $N=400$ we also zoomed into $C$ by dividing it by 2. In the Figure we can see that the combined estimator has better precision everywhere, especially in the regions where the variance of $M$ and $R$ are similar \citep[\emph{cf} Figure 2 of][]{km22}.

\subsection{Asymptotics}

For the asymptotic analysis we employ the following asymptotic expansion of the hypergeometric functions
\begin{align}
N(N-1)z^{2}(1-z)^{N-2}F\left([1, 1, 2-N], [2, 2], -\frac{z}{1-z}\right) & = 1 + \frac{1}{Nz} + \ldots 
\end{align}
so that the asymptotics for the covariance becomes 
\begin{align}
C[R,M] \cdot N = &(q_1 + (2p_Z - 1)C)(1-q_1 - (2p_Z - 1)C) \nonumber \\
& + (q_0 + (2p_Z - 1)C)(1-q_0 - (2p_Z - 1)C) \nonumber \\
& + \frac{(1-p_X)p_Z(q_1 + C - (1-p_Z)D)(1 - q_1 - C + (1-p_Z)D)}{pX}\left(1 +\frac{1}{Np_X}\right)
\nonumber \\
& + \frac{(1-p_X)(1-p_Z)(q_1 - C +p_Z D)(1 - q_1 + C -p_Z D)}{pX}\left(1 +\frac{1}{Np_X}\right)
\nonumber \\
& + \frac{p_Xp_Z(q_0 + C + (1-p_Z)D)(1 - q_0 - C - (1-p_Z)D)}{1-pX}\left(1 +\frac{1}{N(1-p_X)}\right)
\nonumber \\
& + \frac{p_X(1-p_Z)(q_0 - C - p_Z D)(1 - q_0 + C + p_Z D)}{1-pX}\left(1 +\frac{1}{N(1-p_X)}\right)
\nonumber \\
& - 2p_Z(1-p_Z)(4C^2 - D^2) + O(N^{-2})
\end{align}
In the scaling limit $C \sim N^{-\frac{1}{2}}$ and $D \sim N^{-\frac{1}{2}}$ this reduces, rather neatly, to
\begin{align}
\left(C[R,M] -V[R]\right)\cdot N &= - \frac{p_Z(1-p_Z)}{p_X(1-p_X)}\left[2C + \left(2p_X - 1\right)D\right]^2 + O(N^{-\frac{3}{2}}) \\
\left(C[R,M] -V[M]\right)\cdot N &= - \frac{q_1(1-q_1)(1-p_X)}{Np_X^2} - \frac{q_0(1-q_0)p_X}{N(1-p_X)^2} + O(N^{-\frac{3}{2}})
\end{align}

\section{Conclusions}

Both leading-order terms above are non-positive, so in general $C[R, M] < \min(V[R], V[M])$, as also observed in Figure \ref{fig:deltacoplot}. This condition means that $V[P(\alpha^*)] < \min(V[R], V[M])$, so that there exists an estimator with better precision than either $R$ or $M$.

\vskip 0.2in
\bibliography{Vcausal}

\appendix

\newpage

\section{Variances when $X$ is fixed}

Here we repeat the calculations in the setting where node $X$ is fixed by block randomisation, so that the number of entries with $X=0$ is $N^0$ while the number with $X=1$ is $N^1$. We define
\begin{equation}
\xi = \frac{N^1}{N}
\end{equation}
as the fixed fraction of cases where $X=1$ as the analogue of $p_X$ when $X$ was random. The probability tables for $Z$ and $Y$ remain unchanged, but now when sampling $N$ binary vectors, from the DAG in Figure \ref{fig:vstruct} we sample directly from two different multinomials with probabilities
\begin{equation}
\begin{array}{c|c|c|lcc|c|c|l}
X & Z & Y & p & & X & Z & Y & p \\
\cline{1-4}
\cline{6-9}
0 & 0 & 0 & \rho_0 = (1-p_Z)(1-p_{Y,0}) & & 1 & 0 & 0 & \rho_4 = (1-p_Z)(1-p_{Y,2}) \\
0 & 0 & 1 & \rho_1 = (1-p_Z)p_{Y,0} & & 1 & 0 & 1 & \rho_5 = (1-p_Z)p_{Y,2} \\
0 & 1 & 0 & \rho_2 = p_Z(1-p_{Y,1}) & & 1 & 1 & 0 & \rho_6 = p_Z(1-p_{Y,3}) \\
0 & 1 & 1 & \rho_3 = p_Zp_{Y,1} & & 1 & 1 & 1 & \rho_7 = p_Zp_{Y,3}
\end{array}
\end{equation}
where we use $\rho$ to distinguish from the previous probabilities $p$ without fixed numbers in each category of $X$.
If we again represent with $N_i$ the number of sampled binary vectors indexed by $i = 4X+2Z+Y$, then the estimator of $F$ from the raw conditionals is simply
\begin{align}
R = R_1 - R_0\, , & & R_1 = \frac{N_5 + N_7}{N^1}\, , & & R_0 = \frac{N_1 + N_3}{N^0}
\end{align}

Using the marginalisation we would have the following estimator
\begin{align}
M = M_1 - M_0 \, , & &
M_1 = M_{11} + M_{10} \, , & &  M_0 = M_{01} + M_{00}
\end{align}
with the terms separated for later ease
\begin{align}
M_{11} = \frac{N_7}{(N_6+N_7)}\frac{(N_2 + N_3 + N_6 + N_7)}{N} \, , & & M_{01} = \frac{N_3}{(N_2+N_3)}\frac{(N_2 + N_3 + N_6 + N_7)}{N} \nonumber \\
M_{10} = \frac{N_5}{(N_4+N_5)}\frac{(N_0 + N_1 + N_4 + N_5)}{N} \, , & & M_{00} = \frac{N_1}{(N_0+N_1)}\frac{(N_0 + N_1 + N_4 + N_5)}{N}
\end{align}

\subsection{Raw Conditionals}

To compute $E[R]$ we average over a multinomial sample
\begin{equation}
E[R] = \sum \frac{N^0!N^1!}{N_0! \cdots N_7!}\rho_0^{N_0}\cdots \rho_7^{N_7} R
\end{equation}
for which we use that fact that $(\rho_0 + \ldots + \rho_3)^{N^0}(\rho_4 + \ldots + \rho_7)^{N^1}$ generates the probability distribution when we perform a multinomial expansion. To obtain the terms needed for the expectation we define
\begin{align} \label{SNcomp}
S_{N^0, N^1} = & \left[\rho_0+\rho_2 + (\rho_1+\rho_3)v_1\right]^{N^0}\left[\rho_4+\rho_6 + (\rho_5+\rho_7)v_2\right]^{N^1}
\end{align}
with the introduction of two auxilliary generating variables $v_1,v_2$ and whose expansion is
\begin{align} \label{SNmulti}
S_{N^0, N^1} = & \sum \frac{N^0!N^1!}{N_0! \cdots N_7!}\rho_0^{N_0}\cdots \rho_7^{N_7} v_1^{N_1 + N_3} v_2^{N_5 + N_7}
\end{align}

With our generating functions, we can express expectations in terms of differential and integral operators \citep{wilf2005}. For example the operator $v_2 \frac{\partial}{\partial v_2}$ acting on $S_{N^0, N^1}$ will bring down a factor of $(N_5 + N_7)$ as can be seen when we apply it to the expanded form of $S_{N^0, N^1}$ from Equation (\ref{SNmulti}):
\begin{align}
v_2 \frac{\partial}{\partial v_2} S_{N^0, N^1} = & \sum \frac{N^0!N^1!}{N_0! \cdots N_7!}\rho_0^{N_0}\cdots \rho_7^{N_7} (N_5 + N_7) v_1^{N_1 + N_3}v_2^{N_5 + N_7}
\end{align}
Removing the generating variables after applying the operator leads to
\begin{align}
v_2 \frac{\partial}{\partial v_2} S_{N^0, N^1} \Big|_{\substack{v_1=v_2=1}} = & \sum \frac{N^0!N^1!}{N_0! \cdots N_7!}\rho_0^{N_0}\cdots \rho_7^{N_7} (N_5 + N_7) = E[N_5 + N_7]
\end{align}
which is an expectation over the multinomial probability distribution of a binary v-structure with fixed $X$. To actually perform the differentiation we employ the compact form of $S_N$ from Equation (\ref{SNcomp}) to easily obtain the result of $N^1(\rho_5 + \rho_7)$, or
\begin{align}
E[R_1] = \frac{v_2}{N^1}\frac{\partial}{\partial v_2} S_{N^0, N^1} \Bigg|_{\substack{v_1=v_2=1}} = v_2(\rho_5+\rho_7)S_{N^0, N^1 - 1} \, \Bigg|_{\substack{v_1=v_2=1}} = (\rho_5+\rho_7)
\end{align}
Performing the same steps for $R_0$ we obtain
\begin{align}
E[R] = (\rho_5+\rho_7) - (\rho_1+\rho_3)
\end{align}

\sbsbsection{The Variance}
To compute the variance 
\begin{equation}
V[R] = V[R_1] - 2C[R_1, R_0] + V[R_0]
\end{equation}
we first show that the covariance is 0 (as it should be from the sampling setup)
\begin{align}
E[R_1R_0] &= \frac{v_1}{N^0} \frac{\partial}{\partial v_1} \frac{v_2}{N^1}\frac{\partial}{\partial v_2} S_{N^0, N^1} \Bigg|_{\substack{v_1=v_2=1}} = v_1(\rho_1+\rho_3)v_2(\rho_5+\rho_7)S_{N^0 -1, N^1 -1}  \Bigg|_{\substack{v_1=v_2=1}} \nonumber \\
&= E[R_1]E[R_0] 
\end{align}
The last equality follows by comparing to the values of $E[R_1]$ and $E[R_0]$ computed above.

Next we have
\begin{align}
E[R_1^2] &=\frac{v_2}{N^1}\frac{\partial}{\partial v_2}\frac{v_2}{N^1}\frac{\partial}{\partial v_2} S_{N^0, N^1} \Bigg|_{\substack{v_1=v_2=1}} = \frac{v_2}{N^1}\frac{\partial}{\partial v_2} v_2(\rho_5+\rho_7) S_{N^0, N^1 -1}\Bigg|_{\substack{v_1=v_2=1}} \nonumber \\
&= \frac{v_2(\rho_5+\rho_7)}{N^1} S_{N^0, N^1 -1} + \frac{v_2^2(\rho_5+\rho_7)^2(N^1 - 1)}{N^1} S_{N^0, N^1 -2} \Bigg|_{\substack{v_1=v_2=1}} \nonumber \\
&= \frac{(\rho_5+\rho_7)(1 - \rho_5 - \rho_7)}{N^1} + E[R_1]^2
\end{align}
which is just the standard variance of a binomial distribution. Repeating the calculations for $V[R_0]$ we obtain
\begin{align}
V[R] = & \frac{(\rho_5+\rho_7)(1 - \rho_5 - \rho_7)}{N^1} + \frac{(\rho_1+\rho_3)(1 - \rho_1 - \rho_3)}{N^0}
\end{align}

\subsection{Marginalisation}

To compute the expected value $E[M]$  we define
\begin{align}
T_{N^0, N^1} = & \left[\rho_0s_1t_1 + \rho_1s_1t_1u_1 + \rho_2s_2t_2 + \rho_3s_2t_2u_2\right]^{N^0}\left[\rho_4s_1t_3 + \rho_5s_1t_3u_3 + \rho_6s_2t_4 + \rho_7s_2t_4u_4\right]^{N^1}
\end{align}
where we include extra generating variables for all terms in our estimators for which we need the ten generating variables $s_1,s_2,t_1,t_2,t_3,t_4,u_1,u_2,u_3,u_4$. Then
\begin{align} \label{Mexpected}
E[M_{11}] &= \frac{1}{N}\left[\int\mathrm{d}t_4 \frac{s_2u_4}{t_4}\frac{\partial^2}{\partial s_2\partial u_4} \right] T_{N^0, N^1} \, \Bigg|_{\substack{s_1=s_2=1 \cr t_1=t_2=t_3=t_4=1 \cr u_1=u_2=u_3=u_4=1}} \nonumber \\
&= \frac{s_2u_4\rho_7}{N}\left[\frac{t_2(\rho_2+\rho_3u_2)}{(\rho_6 + \rho_7u_4)}N^0T_{N^0 -1, N^1} + t_4 N^1T_{N^0, N^1 -1}\right] \, \Bigg|_{\substack{s_1=s_2=1 \cr t_1=t_2=t_3=t_4=1 \cr u_1=u_2=u_3=u_4=1}} \nonumber \\
&= \frac{\rho_7}{(\rho_6 + \rho_7)}p_Z = \rho_{7}
\end{align}
and similarly for the other terms, leading to
\begin{align}
E[M] = & \rho_7 + \rho_5 - \rho_3 - \rho_1
\end{align}

To compute the variance, we reapply the operators of Equation (\ref{Mexpected}).

\sbsbsection{A Variance}
For example, removing generating variables that do not play a role:
\begin{align}
E[M_{11}^2] \cdot N^2 = & \int\mathrm{d}t_4 \frac{s_2u_4}{t_4}\frac{\partial^2}{\partial s_2\partial u_4} \frac{s_2u_4\rho_7p_Z}{(\rho_6 + \rho_7u_4)}N^{0}T_{N^0 -1, N^1} \, \Bigg|_{\ldots = 1} \nonumber \\
& +  \int\mathrm{d}t_4 \frac{\partial^2}{\partial s_2\partial u_4}s_2u_4\rho_7N^{1}T_{N^0, N^1-1} \, \Bigg|_{\ldots = 1}
\end{align}
For the second line we integrate first to get
\begin{align}
\rho_{7}N + \rho_{7}^2N^1(N - 1) - \frac{\rho_{7}^2}{p_Z}N^{0}
\end{align}
while for the rest of $E[M_{11}^2]$ we first differentiate wrt $s_2$
\begin{align}
s_2\frac{\partial}{\partial s_2}\left[\frac{s_2u_4\rho_7p_Z}{(\rho_6 + \rho_7u_4)}\right] T_{N^0 -1, N^1} \, \Bigg|_{\ldots = 1}
= & \frac{u_4\rho_7p_Z}{(\rho_6 + \rho_7u_4)}\left[T_{N^0 -1, N^1} + (N^{0}-1)p_ZT_{N^0 -2, N^1}\right] \nonumber \\
& {} + t_4N^{1}u_4\rho_7p_ZT_{N^0 -1, N^1 -1} \nonumber \\
\end{align}
For the part with the factor of $t_4$, we again integrate first wrt $t_4$ and then differentiate to obtain
\begin{align}
\rho_{7}N^{0}\left(1 - \frac{\rho_{7}}{p_Z} + \rho_{7}N^{1}\right)
\end{align}
on the rest we apply the operator for $u_4$
\begin{align}
u_4\frac{\partial}{\partial u_4} \ldots \Bigg|_{u_4=1} = &N^{0}\frac{\rho_{6}\rho_{7}}{p_{Z}}T_{N^{0}-1, N^1} + 
N^{0}(N^{0}-1)\rho_{6}\rho_{7}T_{N^{0}-2, N^1} \nonumber \\
& + yN^{0}N^{1}\rho_{7}^2T_{N^{0}-1, N^1 - 1} + yN^{0}N^{1}(N^{0}-1)\rho_{7}^2p_ZT_{N^{0}-2, N^1 - 1}
\end{align}
The linear terms in $t_4$ give the following
\begin{align}
\rho_7^2N^0(N^0 - 1) + \frac{\rho_{7}^2}{p_Z}N^{0} 
\end{align}
while the integrals lead to
\begin{align}
& N^0N^1\rho_6\rho_7(1 + (N^0-1)p_Z)(1-p_Z)^{N^1-1} F\left([1, 1, 1-N^1], [2, 2], -\frac{p_Z}{1-p_Z}\right)
\end{align}
Combining all the terms, subtracting the mean part squared and simplifying slightly we obtain
\begin{align}
& V[M_{11}] \cdot N^2 = N\rho_7(1-\rho_7) + \frac{\rho_6\rho_7}{p_Z}N^{0} \nonumber \\ 
& {} + N^0N^1\rho_6\rho_7(1 + (N^0-1)p_Z)(1-p_Z)^{N^1-1} F\left([1, 1, 1-N^1], [2, 2], -\frac{p_Z}{1-p_Z}\right)
\end{align}

\sbsbsection{Covariances}
For the covariances where separate generating variables are used
\begin{align}
E[M_{11}M_{10}] = \frac{s_2u_4\rho_7}{N^2}\int\mathrm{d}t_3 \frac{s_1u_3}{t_3}\frac{\partial^2}{\partial s_1\partial u_3}\left[\frac{t_2(\rho_2+\rho_3u_2)}{(\rho_6 + \rho_7u_4)}N^0T_{N^0 -1, N^1} + t_4 N^1T_{N^0, N^1 -1}\right]\Bigg|_{\ldots = 1}
\end{align}
the operators act on the $T$ repeating the calculations for the mean with reduced $N$ giving
\begin{align}
C[M_{11},M_{10}] = -\frac{1}{N}E[M_{11}]E[M_{10}] \, , & & C[M_{01},M_{10}] = -\frac{1}{N}E[M_{01}]E[M_{10}] \nonumber \\ 
C[M_{11},M_{00}] = -\frac{1}{N}E[M_{11}]E[M_{00}] \, , & & C[M_{01},M_{00}] = -\frac{1}{N}E[M_{01}]E[M_{00}]
\end{align}
The more complicated cases are where the generating variables reoccur
\begin{align}
E[M_{11}M_{01}] = \frac{u_4\rho_7}{N^2}\int\mathrm{d}t_2 \frac{s_2u_2}{t_2}\frac{\partial^2}{\partial s_2\partial u_2}s_2\left[\frac{t_2(\rho_2+\rho_3u_2)}{(\rho_6 + \rho_7u_4)}N^0T_{N^0 -1, N^1} + t_4 N^1T_{N^0, N^1 -1}\right]\Bigg|_{\ldots = 1}
\end{align}
For the term linear in $t_2$ we first integrate then differentiate wrt $u_2$ while for the other term we first differentiate then integrate to give
\begin{align}
E[M_{11}M_{01}] & = E[M_{11}]E[M_{01}] + \frac{\rho_3\rho_7}{Np_{Z}}\left(1 - p_{Z}\right)
\end{align}
and
\begin{align}
C[M_{11},M_{01}] &= \frac{\rho_3\rho_7(1-p_Z)}{Np_{Z}} \nonumber \\
C[M_{01},M_{00}] &= \frac{\rho_1\rho_5p_Z}{N(1-p_{Z})}
\end{align}

\sbsbsection{The Variance}
The complete variance is
\begin{align}
& V[M] \cdot N^2 = N^0N^1\rho_6\rho_7(1 + (N^0-1)p_Z)(1-p_Z)^{N^1-1} F\left([1, 1, 1-N^1], [2, 2], -\frac{p_Z}{1-p_Z}\right) \nonumber \\ 
& + N^0N^1\rho_4\rho_5(1 + (N^0-1)(1-p_Z))(p_Z)^{N^1-1} F\left([1, 1, 1-N^1], [2, 2], -\frac{1-p_Z}{p_Z}\right) \nonumber \\ 
& + N^0N^1\rho_2\rho_3(1 + (N^1-1)p_Z)(1-p_Z)^{N^0-1} F\left([1, 1, 1-N^0], [2, 2], -\frac{p_Z}{1-p_Z}\right) \nonumber \\ 
& + N^0N^1\rho_0\rho_1(1 + (N^1-1)(1-p_Z))(p_Z)^{N^0-1} F\left([1, 1, 1-N^0], [2, 2], -\frac{p_Z}{1-p_Z}\right) \nonumber \\  
& + N\rho_7(1-\rho_7) + N\rho_5(1-\rho_5) + N\rho_3(1-\rho_3)  + N\rho_1(1-\rho_1)   \nonumber \\
& +  \frac{\rho_6\rho_7}{p_Z}N^{0} +  \frac{\rho_4\rho_5}{1-p_Z}N^{0} +  \frac{\rho_2\rho_3}{p_Z}N^{1} +  \frac{\rho_0\rho_1}{1-p_Z}N^{1} \nonumber \\
& + 2N(\rho_7 - \rho_3)(\rho_1 - \rho_5) -  2N\frac{\rho_3\rho_7(1-p_Z)}{p_{Z}} - 2N \frac{\rho_1\rho_5p_Z}{(1-p_{Z})}
\end{align}

\subsection{The covariance of $R$ and $M$}

For this covariance we also need the two generating variables we had for $R$ and the generating function becomes
\begin{align}
U_{N^0, N^1} = & \left[\rho_0s_1t_1+\rho_2s_2t_2 + (\rho_1s_1t_1u_1+\rho_3s_2t_2u_2)v_1\right]^{N^0} \\ \nonumber
& \left[\rho_4s_1t_3+\rho_6s_2t_4 + (\rho_5s_1t_3u_3+\rho_7s_2t_4u_4)v_2\right]^{N^1}
\end{align}
The calculations proceed as before, for example
\begin{align}
E[M_{11}R_{1}] &= \frac{v_2}{N^1}\frac{\partial}{\partial v_2}\frac{s_2u_4\rho_7v_2}{N}\left[\frac{t_2(\rho_2+\rho_3u_2v_2)}{(\rho_6 + \rho_7u_4v_2)}N^0U_{N^0 -1, N^1} + t_4 N^1U_{N^0, N^1 -1}\right] \, \Bigg|_{\substack{s_1=s_2=1 \cr t_1=t_2=t_3=t_4=1 \cr u_1=u_2=u_3=u_4=1 \cr v_1=v_2=1}}
\end{align}
where we first applied the operators for $M$ as above. With the final differentiation we obtain
\begin{align}
E[M_{11}R_{1}] &= E[M_{11}]E[R_{1}]\left(1-\frac{1}{N}\right) + \frac{\rho_6\rho_7}{p_Z N^1} + \frac{\rho_7^2}{p_Z N} 
\end{align}
Differentiating instead with respect to $v_1$ we end up with
\begin{align}
E[M_{11}R_{0}] &= E[M_{11}]E[R_{0}]\left(1-\frac{1}{N}\right) + \frac{\rho_3\rho_7}{p_Z N} 
\end{align}

Similarly
\begin{align}
C[M_{00},R_{1}] & = -\frac{1}{N}E[M_{00}]E[R_{1}] + \frac{\rho_1\rho_5}{(1-p_Z) N}  \nonumber \\
C[M_{00},R_{0}] & = -\frac{1}{N}E[M_{00}]E[R_{0}] + \frac{\rho_0\rho_1}{(1-p_Z) N^0} + \frac{\rho_1^2}{(1-p_Z) N} \nonumber \\
C[M_{01},R_{1}] & = -\frac{1}{N}E[M_{01}]E[R_{1}] + \frac{\rho_3\rho_7}{p_Z N} \nonumber \\
C[M_{01},R_{0}] & = -\frac{1}{N}E[M_{01}]E[R_{0}] + \frac{\rho_2\rho_3}{p_Z N^0} + \frac{\rho_3^2}{p_Z N}  \nonumber \\
C[M_{10},R_{1}] & = -\frac{1}{N}E[M_{10}]E[R_{1}] + \frac{\rho_4\rho_5}{(1-p_Z) N^1} + \frac{\rho_5^2}{(1-p_Z) N} \nonumber \\
C[M_{10},R_{0}] & = -\frac{1}{N}E[M_{10}]E[R_{0}] + \frac{\rho_1\rho_5}{(1-p_Z) N}  \nonumber \\
\end{align}
so the covariance therefore reduces to
\begin{align}
C[M, R] = & -\frac{(\rho_7 + \rho_5 - \rho_3 - \rho_1)^2}{N} + \frac{(\rho_7 - \rho_3)^2}{p_Z N} + \frac{(\rho_5 - \rho_1)^2}{(1-p_Z)N} \nonumber \\
& + \frac{\rho_6\rho_7}{p_Z N^1} +  \frac{\rho_4\rho_5}{(1-p_Z) N^1} + \frac{\rho_2\rho_3}{p_Z N^0} + \frac{\rho_0\rho_1}{(1-p_Z) N^0}
\end{align}

\subsection{Numerical checks}

Code to evaluate the variance of the two estimators and their covariance through simulation, as well as to evaluate the analytical results for fixed $X$ is also hosted at \url{https://github.com/jackkuipers/Vcausal}.

As an example, for $N=99$ with $N^0 = 66$ and $N^1 = 33$ (corresponding to $\xi = \frac{1}{3}$) and setting $p_Z = \frac{2}{3}, p_{Y,0} = \frac{1}{6}, p_{Y,1} = \frac{1}{2}, p_{Y,2} = \frac{1}{3}, p_{Y,3} = \frac{5}{6}$ we obtained Monte Carlo estimates of the standard deviation of $R$ and $M$ as $0.101662$ and $0.092129$ respectively from 40 million repetitions, and the square root of their covariance as $0.091318$. This agrees with the respective analytical results of $0.101660$ and $0.092132$ for the standard deviations and $0.091321$ for the square root of the covariance.

Since the covariance is lower than either variance, this means that the best combination would have $\alpha^*=0.069409$, for which we would have a minimum variance of $0.092076$. This aligns with the result from the simulation of $0.092073$.

\subsection{Relative difference plots}

\begin{figure}
 \centering
 \begin{tabular}{cc}
 \includegraphics[width=0.45\textwidth]{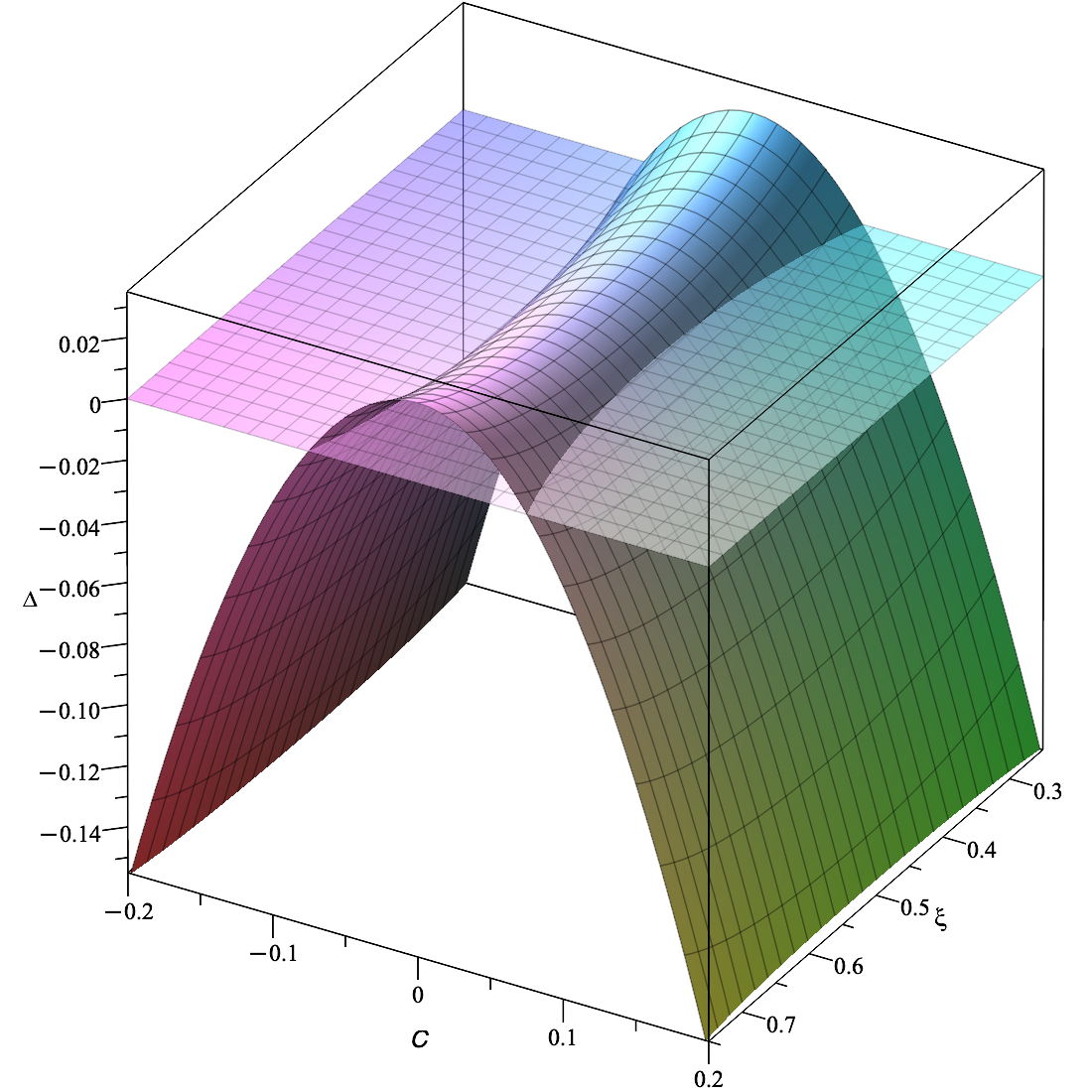} & \includegraphics[width=0.45\textwidth]{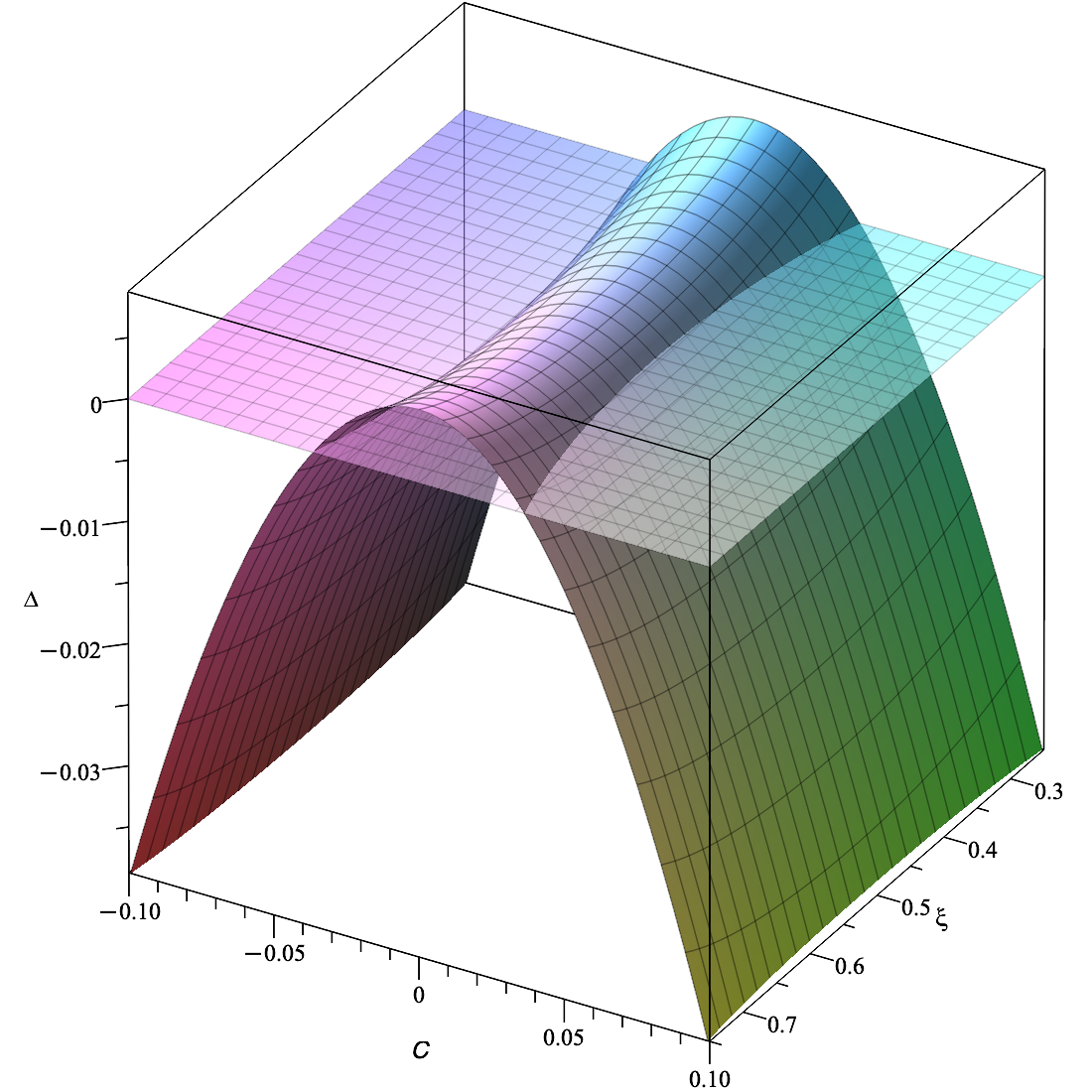} \\
 (a) $N=100$ & (b) $N=400$
 \end{tabular}
\caption{Relative difference in variance of the two estimators. (a) $N=100$, (b) $N=400$.}
\label{fig:deltaplotN}
\end{figure}

With the same parametrisation of $p_Y$ as before, we plot the relative difference in variance of the two estimators $\Delta = \frac{V[M]- V[R]}{V[R]}$ with $\xi$ free and setting $p_Z=\frac{2}{3}$ and $D=0$, along with $q_0=\frac{1}{3}$, $q_1=\frac{2}{3}$ as before. We plot $\Delta$ for $N=100$ and $N=400$ in Figure \ref{fig:deltaplotN}.

\begin{figure}
 \centering
 \begin{tabular}{cc}
 \includegraphics[width=0.45\textwidth]{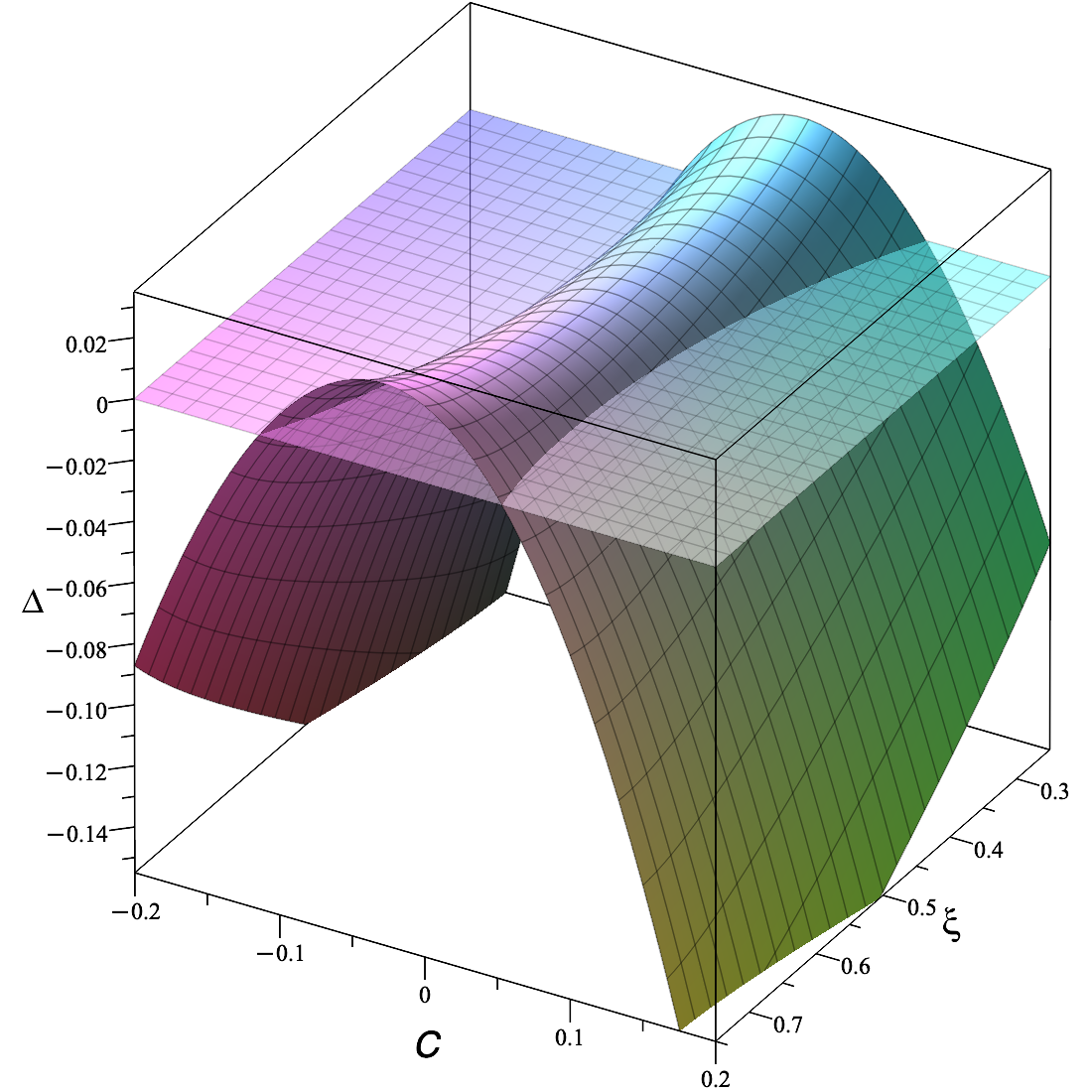} & \includegraphics[width=0.45\textwidth]{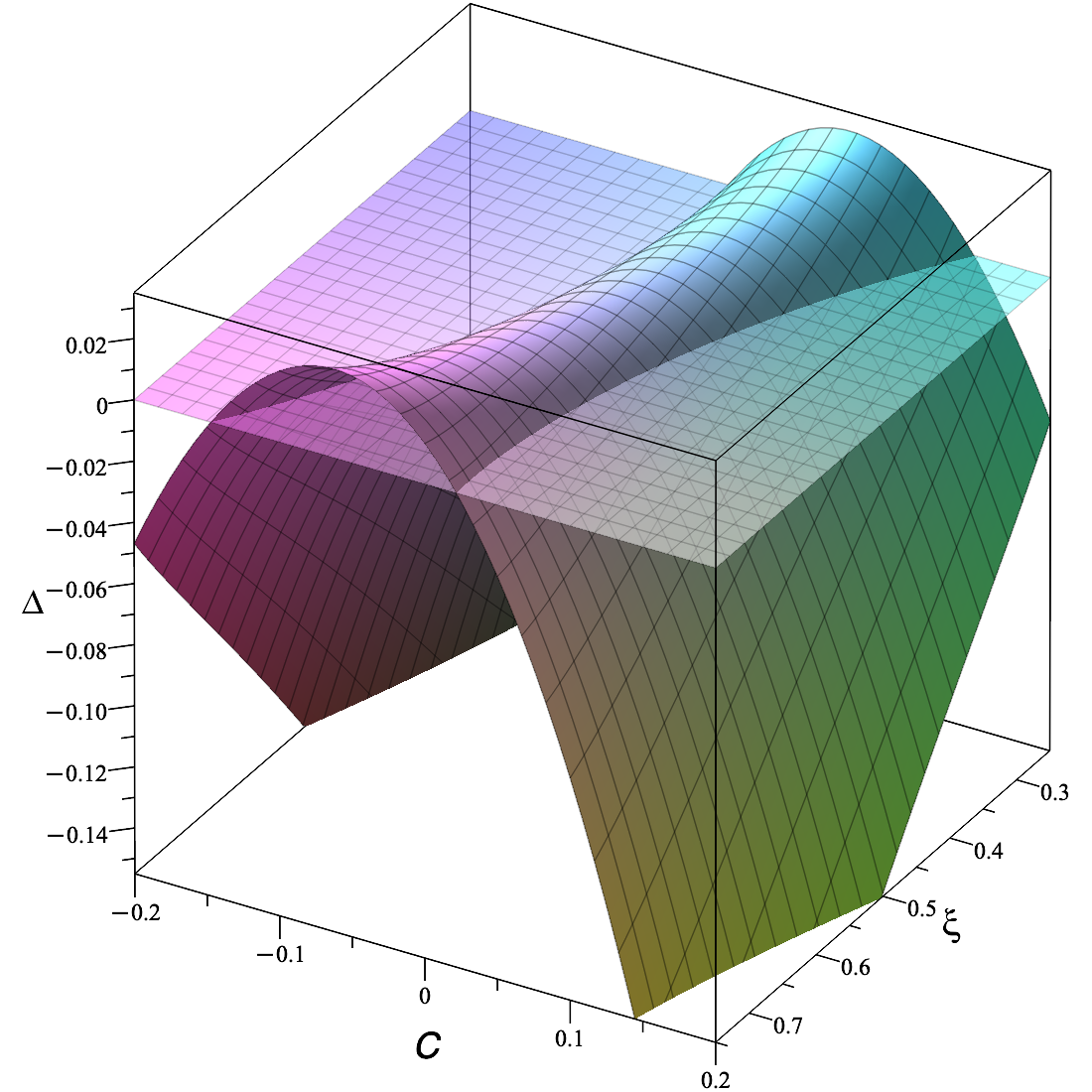} \\
 (a) $D=\frac{1}{8}$ & (b) $D=\frac{1}{4}$
 \end{tabular}
\caption{Relative difference in variance of the two estimators for $N=100$ with interactions. (a) $D=\frac{1}{8}$, (b) $D=\frac{1}{4}$.}
\label{fig:deltaplotD}
\end{figure}

For $N=100$ we also plot the case with an interaction $D\neq 0$ in Figure \ref{fig:deltaplotD}. These plots are highly similar to the case when $X$ is random from \cite{km22}.

\subsection{Asymptotics}

When we consider the asymptotic behaviour in the scaling limit $C \sim N^{-\frac{1}{2}}$ and $D \sim N^{-\frac{1}{2}}$ we obtain that
\begin{align}
\left(V[M] - V[R]\right)\cdot N =& \frac{q_1(1-q_1)(1-\xi)}{\xi^2} + \frac{q_0(1-q_0)\xi}{N(1-\xi)^2} \nonumber \\
&- \frac{p_Z(1-p_Z)}{\xi(1-\xi)}\left[2C + \left(2p_X - 1\right)D\right]^2 + O(N^{-\frac{3}{2}}) \\
\left(C[R,M] - V[R]\right)\cdot N =& - \frac{p_Z(1-p_Z)}{\xi(1-\xi)}\left[2C + \left(2p_X - 1\right)D\right]^2 + O(N^{-\frac{3}{2}}) \\
\left(C[R,M] - V[M]\right)\cdot N =& - \frac{q_1(1-q_1)(1-\xi)}{\xi^2} - \frac{q_0(1-q_0)\xi}{N(1-\xi)^2} + O(N^{-\frac{3}{2}})
\end{align}
which are the same as the results for random $X$ just with a change of $p_X$ for $\xi$.

As for that case covered in the main text, also when $X$ is fixed we can see that the covariance is generally lower than either variance. Again there is therefore generally a linear combination with better precision than both the $R$ and $M$ estimators.

\subsection{Combined estimator plot}

\begin{figure}
 \centering
 \begin{tabular}{cc}
 \includegraphics[width=0.45\textwidth]{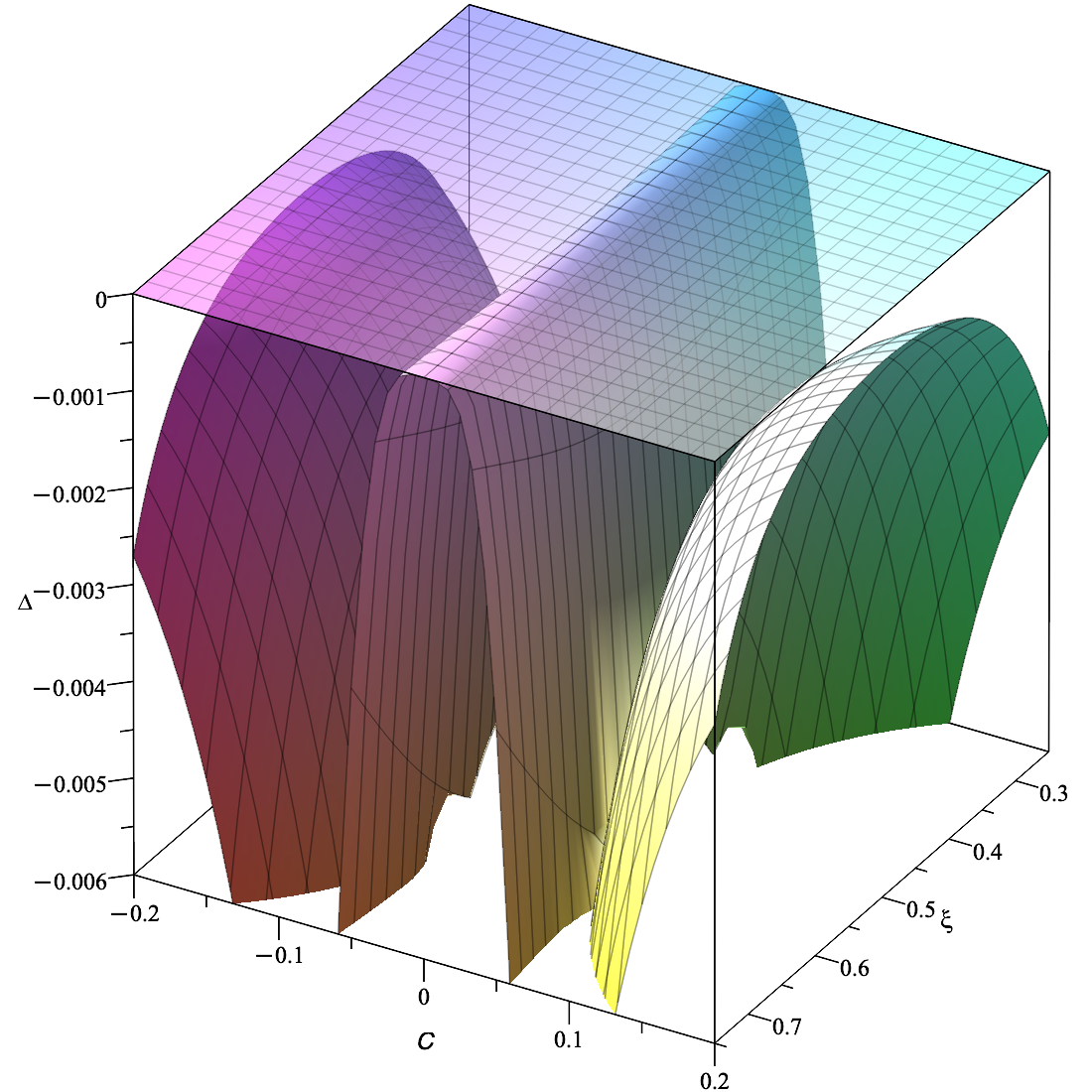} & \includegraphics[width=0.45\textwidth]{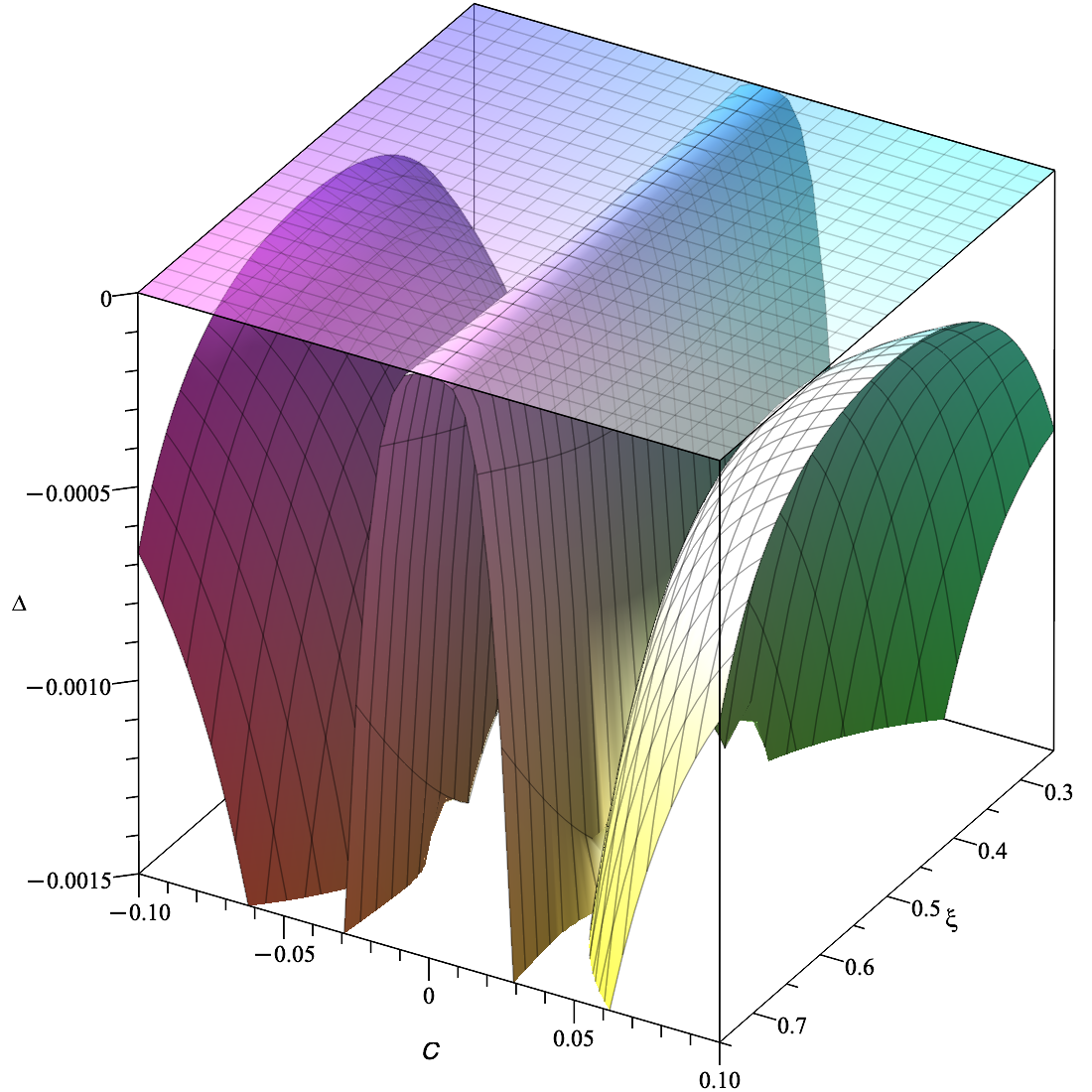} \\
 (a) $N=100$ & (b) $N=400$
 \end{tabular}
\caption{Relative difference in variance of the combined estimator to the better of the other two when $X$ is fixed. This is the analogue of Figure \ref{fig:deltacoplot} in the main text.}
\label{fig:deltacoplotfixed}
\end{figure}

Finally, we can repeat the plot of Figure \ref{fig:deltacoplot} for the case now when $X$ is fixed. With $\xi$ taking the role of $p_X$ and all other parameter values the same, the result in Figure \ref{fig:deltacoplotfixed} is very similar. This is in line with the asymptotic results above also being essentially the same for both cases.

\end{document}